\documentclass[10pt]{article}

\usepackage[english]{babel}
\usepackage{hyperref}
\usepackage[utf8]{inputenc}
\usepackage{amsmath}
\usepackage{graphicx}
\usepackage{amsfonts}
\usepackage[colorinlistoftodos]{todonotes}
\usepackage{fullpage}

\title{Surjectivity of maps induced on matrices by polynomials and entire functions}
\author{Shubhodip Mondal}
\date{March 17, 2015}
\begin{document}
\maketitle

\begin{abstract}
We determine a necessary and sufficient condition for a polynomial
over an algebraically closed field $k$ to induce a surjective map on
matrix algebras $M_n(k)$ for $n \ge 2$. The criterion is given in
terms of algebraic conditions on the polynomial and the proof uses simple linear algebra. Following
that, we formulate and prove a corresponding result for entire
functions as well.
\end{abstract}
\vspace{5 mm}

\begin{description}
\item[Introduction]:

\end{description}
Suppose that $A$ is a complex square matrix. The question of existence of a matrix $B$ such that $B^2 = A$ is a very well-known problem in linear algebra. It turns out that it is not possible to find such a $B$ for every $A$. In this note we extend the above question. Let $f$ be any polynomial over the complex numbers. Is it possible to find a matrix $B$ for every $A$ such that $f(B) = A$? We extend the question even further and ask the same for an entire function $f$ instead of a polynomial. After we prove our results, we can deduce some known results as corollaries e.g., the surjectivity of the exponential map on the
invertible matrices~\cite{ro} and a question (Picard's theorem for matrices) originally asked by Pólya and answered by Szegö, which appears as an exercise in~\cite{ps}.
\vskip 1mm

\noindent We begin with an arbitrary algebraically closed
field $k$. We denote the set of polynomials in one variable over
$k$ by $k[X]$ and the set of $n \times n$ matrices with entries from
$k$ by $\text{M}_{n} (k)$. A polynomial $f \in k[X]$ induces a map
$\text{M}_n(f): \text{M}_n(k) \to \text{M}_n(k)$, where $A \mapsto
\sum a_i A^i \in \text{M}_n(k)$. By abuse of notaton, we denote $\text{M}_n (f) (A)$ as $f(A)$, when $n$ is understood. Let $f' \in k[X]$ denote the
derivative polynomial of $f$ and let $Z(f')$ denote the set of zeros of the polynomial $f'$. If $t \in k$ is such that $f^{-1}(t) \subseteq Z(f')$, we say that $t$ is a critical value of $f$. 
\vskip 5mm

\noindent {\bf Theorem 1.} Let $n \ge 2$ be fixed. $\text{M}_n(f):
\text{M}_n (k) \to \text{M}_n(k)$ is non-surjective iff there exists a $t\in k$ such that $f^{-1}(t) \subseteq Z(f') $, i.e., $f$ has a critical value.\vskip 2mm
\noindent {\bf Note:} The algebraic condition on the polynomial $f$ is independent of $n$. So either $\text{M}_n(f)$ is a surjection for all $n \ge 2$ or a non-surjection for all $n \ge 2$. 
\vskip 2mm
\noindent {\bf Examples:}
\begin{enumerate}
\item Evidently, every polynomial $f \in k[X]$ of degree $1$ induces a
surjection on $M_n(k)$ for any $n \ge 2$.
\item Let $f \in k[X]$ be a quadratic polynomial, i.e., $f = aX^2+bX+c$. If char $k \ne 2$, then $Z(f')$ will be the singleton $\{-b/2a \}$ and the fibre
$f^{-1}(f(-b/2a))= \{-b/2a\}$. If char $k=2$, and $b=0$ then $Z(f')= k$. So applying our result in both of these cases, we see that $M_n(f)$ cannot be a surjection. If char$ k=2$ and $b \ne 0$, then $Z(f')$ is empty, hence $M_n(f)$ is a surjection. 
\item Let char $k = p$. We show that for any $d>2$, there exists a polynomial $f$ of degree $d$ for which $\text{M}_n(f)$ is a surjection. If $p \mid d$, then consider $f(z) = z^d + z$. Since $Z(f')$ is empty, it induces a surjection. If $p$ does not divide $d$ but $p \mid d-1$, then let $f(z) = z^d + z^{d-1}$. Since $Z(f') = 0$ and $\{-1\} \in f^{-1}(0)$, it follows that $f$ induces surjection. The only remaining case is $\gcd(p, d(d-1)) = 1$. Let $f(z) = z^d - dz$. Since $p$ does not divide $d-1$, $f'(z)$ has $d-1$ distinct roots $\{\zeta_1, \ldots,\zeta_{d-1}  \}$. Let there be a $t \in k$ such that $f^{-1}(t) \subseteq Z(f')$. Then for some $1 \le r \le d-1$, $\zeta_r$ is a root of the polynomial $h(z) = f(z) - t$. If the multiplicity of $\zeta_r$ in $h(z)$ is at least 3, the multiplicity of $\zeta_r$ in $h'(z) = f'(z)$ is at least 2. But this is a contradiction since $f'(z)$ has $d-1$ distinct roots. So the multiplicity of $\zeta_r$ in $h(z)$ is at most 2. Since $\deg(h(z)) = d > 2$, there exists another root of $h(z)$. So there exists some  $1 \le r \ne s \le d-1$ such that $\zeta_s$ is a root of $h(z)$. But then $t= f(\zeta_r) = f(\zeta_s)$, which imples $(1-d) \zeta_r = (1-d) \zeta_s$, contradicting the fact that $\zeta_r \ne \zeta_s$. Hence there cannot be any $t \in k$ such that $f^{-1}(t) \subseteq Z(f')$. Therefore $f$ induces a surjection. 

\end{enumerate}

\noindent We recall some definitions before going into the proof of
Theorem 1. Recall that a  Jordan matrix $J$ is a block diagonal
matrix
$$J= \begin{pmatrix}
J_1 & 0 & \cdots & 0 \\ 0 & J_2 & \cdots &  0 \\
\vdots & \vdots & \ddots & \vdots \\
0 & 0 & \cdots & J_p
\end{pmatrix}$$

\noindent where the Jordan blocks $J_i$'s are square matrices of the
form
$$ \begin{pmatrix}
\lambda_i & 1       & 0       & \cdots  & 0 \\
0       & \lambda_i & 1       & \cdots  & 0 \\
\vdots  & \vdots  & \vdots& \ddots  & \vdots \\
0       & 0       & 0        & \lambda_i & 1       \\
0       & 0       & 0       & 0       & \lambda_i \end{pmatrix} $$

\noindent The following facts are quite well-known and their proofs of can be found in~\cite{hk} for example.\\

\begin{enumerate}
\item If $k$ is algebraically closed, any matrix $A \in \text{M}_n(k)$ is similar to a Jordan matrix $J$. The matrix $J$ is said to be the Jordan normal form of $A$.

\item The Jordan normal form of a matrix is unique up to a
permutation of the diagonal blocks. \vspace{2mm}

\item $\lambda_i$'s are eigenvalues of $A$ and number of Jordan blocks
corresponding to the eigenvalue $\lambda_i$ in the Jordan normal
form of $A$ is the dimension of the Kernel of $(A - \lambda_i \cdot
I )$. This is a direct consequence of the rank-nulity theorem. \vskip
3mm
\end{enumerate}

\noindent Now we prove two lemmas which will be used to prove the
theorem.\\

\noindent {\bf Lemma 1.} Let $U$ be a Jordan block with $\lambda$ as
its eigenvalue. Let $p(X) = \sum_{m=0}^{n} a_m X^m \in k[X]$. Then, $p(U)_{ij} = 0$ for $i > j$ and $$p(U)_{ij} = \sum_{m=0}^{n} a_m
\binom{m}{j-i} \lambda ^{m - (j-i)} $$
otherwise. (Here we follow the convention that $\binom{p}{q} = 0$ for $q >p$.) \\

\noindent {\bf Proof.} Observe that $U = \lambda I + N$ where $N_{ij} = 1$ if $j- i = 1$,
and $N_{ij} = 0$ otherwise. By linearity, it suffices to prove the lemma
in the case of $p(X) = X^k$, which is done below. 
\begin{align*}
U^k _{ij} &= (\lambda I + N)^k _{ij}\\
&= \sum_{r=0}^{k}  \binom{k}{r}    \lambda I^r N^{k-r}_{ij}\\
&= \sum_{r=0}^{k} \binom{k} {r} \lambda ^r N^{k-r}_{ij}
\end{align*}
Now, $N^{k-r} _{ij}$ is nonzero, only when $j-i = k-r$, or $r = k -
(j-i)$.\\
Hence the sum equals,
$$ \binom{k} {k - (j-i)} \lambda ^{k- (j-i)} = \binom{k}{j-i}
\lambda ^{k - (j-i)} $$
as asserted. \\

\noindent {\bf Note:}  Irrespective of the characteristic of $k$, if
$j-i = 1$, we have $p(U)_{ij}=p'(\lambda)$. If char $ k = 0$, then $p(U)_{ij} = \frac{p^{(j-i)}(\lambda)}{(j-i)!}$, for all $1 \le
i \le j \le n$ \vskip 3mm.

\noindent {\bf Lemma 2.} Let $n \ge 2$. If $U \in \text{M}_n(k)$ is
a Jordan block with $\lambda$ as its eigenvalue and $p(X) \in k[X]$, then the Jordan
normal form of $p(U)$ has at least two Jordan blocks if and only if
$p'(\lambda) = 0$.\\

\noindent {\bf Proof.} To prove this lemma, we use the third fact noted earlier. By
Lemma 1, $p(\lambda)$ is the only eigenvalue of $p(U)$ and number of
blocks is equal to the dimension of $\text{Ker} ( P(U) -
p(\lambda)\cdot I )$. So if $p' (\lambda) = 0$, then the first two
columns of $P(U) - p(\lambda)\cdot I $ are zero by Lemma 1. Hence the
rank is at most $n-2$. By the rank-nullity theorem, $ \text{Ker} (
P(U) - p(\lambda)\cdot I ) \ge 2$.\\
Conversely, if $p '(\lambda) \ne 0$, the matrix $P(U) -
p(\lambda)\cdot I$ has $n-1$ linearly independent columns. Indeed, by
Lemma 1, the first column $c_1$ of $P(U) - p(\lambda)\cdot I $ is
zero and the $k$-th column (for $2 \le k \le n$) is
$$c_k =  \left( \frac{p^{(k-1)}(\lambda)} {(k-1)!} , \ldots, p'(\lambda), 0 , \ldots, 0 \right)^t.$$
Since $p'(\lambda) \ne 0$, all the $n-1$ vectors $c_k$, for $2 \le k
\le n$ are linearly independent. So $\dim \text{Ker} \left(P(U) -
p(\lambda)\cdot I \right) = 1$, and $p(U)$ has only one Jordan block. \vskip 5mm

\noindent {\bf Proof of Theorem 1.} Let $t$ be a critical value of $f$. We take $Y \in M_n(k)$ such that $Y$ is a Jordan block with $t$ as eigenvalue and show that $Y$ is not in the image of $ M_n(f)$. Assume that $f(X') = Y$. Let $X$ be the Jordan normal form of $X'$. So $f(X)$ is similar to $Y$ (using $f(PAP^{-1}) = P f(A) P^{-1} $) . 
Since $Y$ is a Jordan block, $X$ also has to be a Jordan
block. Otherwise Jordan normal form of $f(X)$ will not be a single Jordan
block and hence cannot be similar to $Y$. Let $u$ be the eigenvalue of $X$. Hence $f(u)$ is the only
eigenvalue of $f(X)$ (by Lemma 1). Since $f(X)$ is similar to $Y$,
there eigenvalues have to be the same. So $f(u) = t$ $
\implies u \in f^{-1} (t) \subset Z(f')$. Hence $f'(u) = 0$. But
then by Lemma 2, Jordan normal form of $f(X)$ has more than one
Jordan block. Therefore, it cannot be similar to a single Jordan block
$Y$. \vskip 2mm

\noindent Before proving the converse, we prove two more lemmas which will be useful later. \vskip 2mm

\noindent {\bf Lemma 3.} Let $r \ge 1$ and $Y \in \text{M}_r(k)$ be a Jordan
block with eigenvalue $\lambda$ such that $\lambda$ is not a critical value of $f$. Then there exists a Jordan block $X$ in $M_r(k)$ such that Jordan form of
$f(X) = Y $.\\

\noindent {\bf Proof.} For $r=1$, the result is clear, since $k$ is algebraically closed.
So we assume that $r \ge 2$ in what follows. Since $f^{-1} (\lambda) \not\subset Z(f')$, there exists $u \in k$ such
that $f(u) = \lambda$ and $f'(u) \ne 0$. Let $X \in \text{M}_r (k)$ be the Jordan block with eigenvalue $u$. By Lemma 2, Jordan form of $f(X)$
is a Jordan block of order $r$ and has eigenvalue $\lambda$, so it
has to be equal to $Y$.\\
 
\noindent {\bf Lemma 4.} Let $B \in M_n(k)$ such that the eigenvalues of $B$ are not critical values of $f$. Then $B$ lies in the image of $M_n(f)$.\\

\noindent {\bf Proof.} Let $Y$ be the Jordan normal form of $B$. Then $Y = \text{diag} (Y_1, \ldots, Y_p)$ is a Jordan matrix, where
$Y_i$'s are Jordan blocks. Using the hypothesis on eigenvalues of $B$ and Lemma 3, there exists a Jordan
block $X_i$ such that Jordan form of $f(X_i) = Y_i$. So Jordan form
of $\text{diag} (f(X_1), \ldots, f(X_p) ) = Y$. Now consider the matrix $X =\text{diag}(X_1, \ldots,X_p) $. Since
$$f(X) = f(\text{diag}(X_1, \ldots,X_p)) = \text{diag} (f(X_1), \ldots, f(X_p) ),$$
the Jordan form of $f(X)$ is equal to $Y$. In particular $f(X)$ is similar to $Y$ and consequently to $B$, which implies that $B$ is in the image of $M_n(f)$ (using $f(PAP^{-1}) = P f(A) P^{-1}$ ). \vskip 2mm
Therefore, if $f$ has no critical values, Lemma 4 implies that $M_n(f)$ is a surjection. This proves the converse and finishes the proof of Theorem 1.\vskip 1mm
We have established the following in the first part of the proof of Theorem 1, which we note down explicitely as a lemma for later use: \vskip 1mm
\noindent {\bf Lemma 5.}  If $t$ is a critical value of $f$ and $Y \in M_n(k)$ is a Joran block with eigenvalue $t$, then $Y$ does not belong to the image of $M_n(f)$.
 
\vskip 2mm

\noindent {\bf Matrices of entire functions:} \vskip 1mm

\noindent A function $f: \mathbb C \to \mathbb C$ is said to be entire if it holomorphic
on the whole of $ \mathbb {C}$. Given an entire function, we write
it as a power series around zero : $f(z) = \sum_{n =0}^{\infty} a_n z^n$. For an entire function $f$,  $ f(A) := \sum_{n=0}^{\infty} a_n A^n  $ is a well-defined matrix. The reader may refer to~\cite{mf} for its proof. So sending $A$ to $f(A)$ gives a map $\text{M}_n(f) : \text{M}_n(\mathbb C) \to \text{M}_n (\mathbb C)$. Let $f'$ be the derivative of $f$ and let $Z(f')$ denote the zeros of the $f'$. If $t \in \mathbb C$ is such that $f^{-1}(t) \subseteq Z(f')$, we say that $t$ is a critical value of $f$. Note that if $t \not\in f(\mathbb C) $, $t$ is a critical value of $f$.

 We define $D_n(f)$ to be the set of all matrices in $M_n(\mathbb C)$ whose eigenvalues lie in $f( \mathbb{C})$. Lemma 1 easily extends to the case of entire functions. In particular, if $\lambda_1, \ldots, \lambda_n$ are eigenvalues of $A \in \text{M}_n (\mathbb C)$, then $f(\lambda_1) ,\ldots,f(\lambda_n) $ are the eigenvalues of $f(A)$. So $\text{Im} ( M_n(f)) \subseteq D(f)$. Since $ f(P A P^{-1}) = P  f (A) P^{-1}$ and $D_n(f)$ is closed under conjugation, the discussion for the case of polynomials applies mutatis mutandis to entire functions and yields the corresponding versions of all our lemmas. We also obtain the following theorem: \vskip 2mm

\noindent {\bf Theorem 2.} Let $n\ge2$ be a fixed. Then $\text{Im} ( M_n(f)) \ne D_n(f)$ iff
there exists a $t \in f( \mathbb C)$ such that $f^{-1}(t) \subseteq
Z(f')$.  \vskip 2mm

 \vspace{1 mm}
\noindent {\bf Note:} By Picard's little theorem~\cite{re}, if $f$ is a non-zero entire function then at most one complex number does not belong to the image of $f$. So $D_n(f)$ is either $\text{M}_n(\mathbb{C})$ or the set of all matrices whose eigenvalues does not equal to $p_f$ for some fixed complex number $p_f$.

\vskip 2 mm
\noindent {\bf Examples.}

\begin{enumerate}
\item Let $f(z) = \text{exp}(z)$. Then $f(\mathbb C) = \mathbb C ^{*}$. Hence $D_n(f) = \text{GL}_n(\mathbb C)$, the set of $n\times n$ invertible matrices over $\mathbb C$. Since $Z( f') = Z(f) = \emptyset$, the matrix exponential is a surjective map
from $M_n(\mathbb C) \to \text{GL}_n (\mathbb C)$.

\item $\text{sin}z$ and $\text{cos}z$ are surjective entire functions. But the maps they induce from $\text{M}_n (\mathbb C) \to \text{M}_n (\mathbb C)$ are not surjective for $n \ge 2$. Indeed, we have $\sin ^2 z + \cos ^2 z = 1$. So $\sin ^{-1} ( \left \{ \pm 1 \right \}) \subseteq Z ( \cos z)$ and $\cos^{-1}( \left \{ \pm 1 \right \})  \subseteq Z(\sin z)$. Hence by Lemma 5 (which remains true for entire functions, as noted in the previous discussion), Jordan blocks with eigenvalue $1$ and $-1$ are not in the image of $\text{M}_n(\text{sin}(z))$ or $\text{M}_n(\text{cos}(z))$. 
\item {\bf Picard's theorem for matrices  ~\cite{ps} :}
Let $C_f$ denote the set of all critical values of an entire function $f$. As a corrollary to the Second Fundamental Theorem of Nevanlinna, one obtains that $|C_f| \le 2$. The reader may refer to~\cite{nv} for an exposition on Nevanlinna theory, which also contains the stated corrollary. Now by using Lemma 4 and Lemma 5 (both of them extend to the case of entire functions, as already noted) we obtain that an entire function $f$ has at most $2$ "exceptional values" in the following special sense: $A \in M_n(\mathbb C)$ lies in the image of $M_n(f)$ if none of the eignevlaues of $A$ coincides with an exceptional value of $f$. On the other hand, there are certain matrices with eigenvalues consisiting of exceptional values not belonging to the image of $M_n(f)$. 

\end{enumerate}
\vskip 2mm
\noindent {\bf Acknowledgements:} I would like to thank Prof. S. Inamdar and
Prof. B. Sury at Indian Statistical Institute for going through the
proof and valuable suggestions.

\end{document}